\documentclass[12pt]{amsart}

\usepackage{latexsym}
\usepackage{a4}
\usepackage{amssymb}
\usepackage{amsmath}
\usepackage[english]{babel}
\usepackage{tikz}
\usepackage{bm}
\usepackage{hyperref}

\newtheorem{thm}{Theorem}[section]
\newtheorem{pro}[thm]{Proposition}
\newtheorem{lem}[thm]{Lemma}

\newtheorem{cor}[thm]{Corollary}
\theoremstyle{definition}

\newtheorem{de}[thm]{Definition}

\numberwithin{equation}{section}

\newcommand{\N}{\mathbb{N}}

\newcommand{\Con}{\operatorname{Con}}

\newcommand{\Pol}{\operatorname{Pol}}

\newcommand{\bS}{\mathbf{S}}
\newcommand{\bA}{\mathbf{A}}
\newcommand{\bB}{\mathbf{B}}

\newcommand{\oa}{\mathbf{a}}
\newcommand{\ob}{\mathbf{b}}
\newcommand{\oc}{\mathbf{c}}
\newcommand{\od}{\mathbf{d}}

\newcommand{\ox}{\mathbf{x}}

\newcommand{\oy}{\mathbf{y}}

\newcommand{\tup}[3]{(#1_{#2},\dots,#1_{#3})}

\title{Supernilpotent Semirings}

\author{Neboj\v{s}a Mudrinski}
\address{Neboj\v{s}a Mudrinski,
Department of Mathematics and Informatics, Faculty of Sciences,
University of Novi Sad, 21000 Novi Sad, Serbia}
\email{\tt
nmudrinski@dmi.uns.ac.rs}

\author{Milica \v Sobot}
\address{Milica \v Sobot,
Universit\"at Bonn, Bonn, Germany} \email{\tt
mica.sobot@gmail.com}

\thanks{Supported by the Ministry of Science, Education and Technological Development of the Republic of Serbia (Grant No. 451-03-68/2022-14/200125)}
\date{\today}

\begin{document}
\bibliographystyle{amsalpha}

\maketitle

\begin{abstract}
We prove that supernilpotent and nilpotent semirings with absorbing zero are the same and provide a necessary and sufficient condition for supernilpotency (nilpotency).
\end{abstract}

\maketitle

\section{Motivation}

In general, an algebra $\bA$ is supernilpotent if there is an $n\in\N$ such that the commutator $[,]$ of the length $n+1$ satisfies $[\,\underbrace{1,\dotsc,1}_{n+1}\,]_{\bA}=0$, where $0$ denotes the smallest (the equality relation) and $1$ is the largest (the full relation) in the congruence lattice $\Con(\bA)$ of the algebra $\bA$. The commutator of any finite length has been introduced by Bulatov in 2001, see Definition \ref{DefnHC}. It is a finitary operation on the congruence lattice of the algebra that extends the term condition commutator of length two, from \cite{FM:CTCMV,McMc}. The binary commutator generalizes the commutator from the group theory. Using the term condition commutator in arbitrary algebra, one can define the properties that exist in group theory that are defined by the group theoretic commutator, such as abelian, nilpotent or solvable, see Definition \ref{DefNilpSolv}. A natural problem is to study what does this general concept of properties defined by the commutator carry back to concrete algebraic structures. In congruence permutable varieties abelian means that the algebra is polynomially equivalent to a module over a ring. A characterization of abelian semigroups has been obtained by Warne in \cite{Warne} and \cite{Warne2}. Solvable and nilpotent regular semigroups have been described in \cite{RM:CCSS}. Supernilpotent orthodox semigroups has been characterized in \cite{RM:ss}. All the properties defined by the commutator have been studied for semigroups with zero in \cite{RM:cswz}.

The next problem that naturally arises is how are the properties defined by the commutator connected? Abelian is a special case of both nilpotency and supernilpotency by definition. Using basic properties of the commutator, see \cite{FM:CTCMV} one can observe that abelian implies nilpotency and nilpotency implies solvability. How is supernilpotency connected with the others? In Taylor algebras supernilpotency
implies nilpotency by Moorhead, see \cite{Moorhead:Taylor}, as well
as in all finite algebras by Kearnes and Szendrei, see \cite{KS:ISSN}. However, by Moore and
Moorhead it is not always true, namely there exist a supernilpotent
infinite algebra that is not nilpotent, see \cite{MM:SNNIN}. 

The aim of this note is to study the two mentioned problems in semirings. For the basic notions and notations in semirings, we follow \cite{HW:BookSemirings}. We are going to see in Proposition
\ref{1} that semirings that have neutral elements for both binary
operations never are abelian, nilpotent or solvable.
Therefore, we omit the neutral element for the second binary
operation in our definition of semiring. Here, a \emph{semiring} is an algebra $\bS=(S,+,\cdot,o)$ such that $(S,+,o)$ is a
commutative monoid, $(S,\cdot)$ is a semigroup with zero $o$ and
$x(y+z)=xy+xz$ and $(x+y)z=xz+yz$ for all $x,y,z\in S$. In \cite{HW:BookSemirings} it is called a semiring with (multiplicatively) \emph{absorbing zero} and in \cite{G:BookSemirings} it is called a hemiring (for example $(2\N_0,+,\cdot,0)$). In such semirings we give the connection between nilpotency and supernilpotency as well as the necessary and sufficient condition for the both of the properties. This is the main result of the note, formulated as follows.
\begin{thm}\label{teoremanilpotent}
Let $n\in\N$ and let $\bS=(S,+,\cdot,o)$ be a semiring. Then the following are equivalent:
\begin{enumerate}
\item $\bS$ is $n$-nilpotent;
\item $\bS$ is additively cancellative and $S^{n+1}=\{o\}$;
\item $\bS$ is  $n$-supernilpotent.
\end{enumerate}
\end{thm}

A semiring $(S,+,\cdot,o)$ is \emph{additively cancellative} if the commutative monoid $(S,+,o)$ is cancellative, see \cite{HW:BookSemirings}. As a special case of Theorem \ref{teoremanilpotent}, we have that a
semiring is abelian if and only if it is additively cancellative and
the product of two elements is always zero, see Proposition
\ref{teoremabelian}. Also, we provide a complete characterization for solvability in additively cancellative semirings, see Proposition \ref{teoremasolvable}. Note that in an additively cancellative semiring $\bS=(S,+,\cdot,o)$, $o$ is always  absorbing zero, see \cite[Fact 2.12. a), p.13]{HW:BookSemirings}.

\section{Ideals and Commutators}

In this section we will first review some definitions and properties about the term condition of algebras in general. For the basic notions such as congruences and polynomials and notations we refer to \cite{FM:CTCMV,McMc}. For an algebra $\bA$, $n\in\N$, tuples $\oa_1,\dots,\oa_n,\ob_1,\dots,\ob_n$ from $A$, tuples of variables $\ox_1,\dots,\ox_n$ such that $\oa_i,\ob_i$ and $\ox_i$ have the same length for all $i\in\{1,\dots,n\}$, polynomial $p\tup{\ox}{1}{n}$ and a valuation $v:\{(\ox_1,\dots,\ox_n)\}\to\{\oa_1,\ob_1\}\times\dots\times\{\oa_n,\ob_n\}$ we denote $p^{\bA}(v(\ox_1),\dots,v(\ox_n))$ by $p^v\tup{\ox}{1}{n}$.

\begin{de}(cf.\cite{Bu:OTNO})\label{DefnC}
Let $\bA$ be an algebra. Let $k\geq 1$ and let $\alpha_1,\dotsc,\alpha_k,\beta$ and $\delta$ be congruences in $\Con(\bA)$. Then we say that $\alpha_1,\dotsc,\alpha_k$ \textit{centralize} $\beta$ \textit{modulo} $\delta$, and we write $C(\alpha_1,\dotsc,\alpha_k,\beta;\delta)_{\bA}$, if for every polynomial $p\in\Pol(\bA)$, of arity $n_1+\dotsc+n_k+m$, where $n_1,\dotsc,n_k,m\in\N$, and every $\oa_1,\ob_1\in A^{n_1}$, $\dotsc$, $\oa_{k},\ob_{k}\in A^{n_{k}}$ and $\oc,\od\in A^{m}$ such that
\begin{enumerate}
\item[(1)] $\oa_j \;\alpha_j\; \ob_j$, for $j=1,\dotsc,k$;
\item[(2)] $\oc \;\beta\; \od$;
\item[(3)] $p^v(\ox_1,\dotsc,\ox_{k},\oc)\;\delta\;p^v(\ox_1,\dotsc,\ox_{k},\od)$ for all $v:\{(\ox_1,\dots,\ox_n)\}\to\{\oa_1,\ob_1\}\times\dots\times\{\oa_n,\ob_n\}$ such that $v\neq w$, where $w(\ox_i)=\ob_i$ for all $i\in\{1,\dots,k\}$;
 \end{enumerate}
we have $p(\ob_1,\dotsc,\ob_{k},\oc) \; \delta \; p(\ob_1,\dotsc,\ob_{k},\od)$.
\end{de}

\begin{de}(\cite{Bu:OTNO})\label{DefnHC}
Let $\bA$ be an algebra. Let $n\geq 2$ and let $\alpha_1,\dotsc,\alpha_n\in\Con(\bA)$. The \textit{$n$-ary commutator of $\alpha_1,\dotsc,\alpha_n$} is the smallest congruence $\delta$ on $\bA$ such that $C(\alpha_1,\dotsc,\alpha_n;\delta)_{\bA}$. We abbreviate it by $[\alpha_1,\dots,\alpha_n]_{\bA}$.
\end{de}

We omit the algebra in the index in the notation of centralizing relation $C(\alpha_1,\dotsc,\alpha_k,\beta;\delta)_{\bA}$ and the commutator $[\alpha_1,\dots,\alpha_n]_{\bA}$ whenever it is clear from the context. Just as an immediate corollary of the previous definitions we have the following statement.

\begin{pro}\label{centralizingreduct}
Let $\bA=(A;F)$ and $\bB=(A;F')$ be two algebras on the same domain $A$ such that for the sets of fundamental operations $F$ and $F'$ we have $F\subseteq F'$. Then for $n\geq 2$ and congruences $\alpha_1,\dots,\alpha_n$ and $\beta$ of $\bB$ (that are also congruences of $\bA$) we have:
\begin{enumerate}
\item $C(\alpha_1,\dots,\alpha_n;\beta)_{\bB}\Rightarrow C(\alpha_1,\dots,\alpha_n;\beta)_{\bA}$;
\item $[\alpha_1,\dots,\alpha_n]_{\bA}\leq[\alpha_1,\dots,\alpha_n]_{\bB}$.
\end{enumerate}
\end{pro}

\begin{proof}
We observe that every polynomial of $\bA$ is also a polynomial of $\bB$. Then we use Definitions \ref{DefnC} and \ref{DefnHC}. 
\end{proof}

In the binary case the commutator can be defined in an equivalent and
more convenient way.

\begin{pro}\label{DefnC(a,a,a,b;d)}
Let $\bA$ be an algebra. Let $\alpha,\beta\in\Con(\bA)$. Then we say
that $\alpha$ \textit{centralizes} $\beta$ \textit{modulo} $\delta$,
and we write $C(\alpha,\beta;\delta)$, if for every polynomial
$p\in\Pol(\bA)$, of arity $n+1$, where $n\in\N$, and every $a,b\in
A$ and $\oc,\od\in A^{n}$ such that $p(a,\oc)\;\delta\; p(a,\od)$ we
have $p(b,\oc) \; \delta \; p(b,\od)$.
\end{pro}

\begin{proof}
See \cite[Exercise 4.156.2]{McMc} and \cite{EA:TPFoCA}.
\end{proof}

\begin{de}\textup{(\cite{FM:CTCMV}, Definition 6.1)}\label{DefNilpSolv}
Let $\bA$ be an algebra, and let $\rho,\sigma$ be congruences on $\bA$. We define the series $(\rho,\sigma]^{(k)}$, $k\in\N$, as follows: $(\rho,\sigma]^{(1)} = [\rho,\sigma]$, and $(\rho,\sigma]^{(k+1)} = [\rho,(\rho,\sigma]^{(k)}]$, for $k\in\N$. Similarly, we define the series $[\rho]^{(k)}, k\in\N$ by $[\rho]^{(1)}=[\rho,\rho]$ and $[\rho]^{(k+1)}=[[\rho]^{(k)},[\rho]^{(k)}]$, for $k\in\N$.

Algebra $\bA$ is \textit{$n$-nilpotent}, $n\in\N$, if $(1,1]^{(n)}=0$. Similarly, an algebra $\bA$ is $n$-solvable, $n\in\N$, if $[1]^{(n)}=0$.
\end{de}

Note that $1$-nilpotency and $1$-solvability is the same property
and we call it \emph{abelian}. 

\begin{pro}\label{forreduct}
Let $\bA=(A;F)$ and $\bB=(A;F')$ be two algebras on the same domain $A$ such that for the sets of fundamental operations $F$ and $F'$ we have $F\subseteq F'$. Then for all $n\in\N$ we have
\begin{enumerate}
\item if $\bB$ is $n$-nilpotent then $\bA$ is $n$-nilpotent;
\item if $\bB$ is $n$-supernilpotent then $\bA$ is $n$-supernilpotent.
\end{enumerate}
\end{pro}

\begin{proof}
(1) One can prove by induction on $k$ that $(1,1]^{(k)}_{\bA}\leq(1,1]^{(k)}_{\bB}$ for all $k\in\N$ using Proposition \ref{centralizingreduct}, Definition \ref{DefNilpSolv} and monotonicity of the commutator.
(2) If $\bB$ is $n$-supernilpotent then we have $C(\,\underbrace{1,\dots,1}_{n+1};0)_{\bB}$ and therefore $C(\,\underbrace{1,\dots,1}_{n+1};0)_{\bA}$ by Proposition \ref{centralizingreduct}.
\end{proof}

\begin{pro}\label{1}
If there is an $e\in S$ in a semiring $\bS=(S,+,\cdot,o)$ such that
$(S,\cdot,e)$ is a monoid then for all $n\geq 2$ we have
$[\,\underbrace{1,\dots,1}_n\,]=1$.
\end{pro}

\begin{proof}
We prove $C(\,\underbrace{1,\dots,1}_n;\theta)\Rightarrow
\theta=1$. We assume $C(\,\underbrace{1,\dots,1}_n;\theta)$
and we take the polynomial $p(x_1,\dots,x_n)=x_1\cdot\ldots\cdot x_n$.
Then we have
$p(\alpha_1,\dots,\alpha_{n-1},o)=o=p(\alpha_1,\dots,\alpha_{n-1},e)$
for all
$(\alpha_1,\dots,\alpha_{n-1})\in\{o,e\}\times\ldots\times\{o,e\}\backslash\{(e,\dots,e)\}$
and therefore
$$
p(\alpha_1,\dots,\alpha_{n-1},o)\,\theta\,
p(\alpha_1,\dots,\alpha_{n-1},e)
$$
for all
$(\alpha_1,\dots,\alpha_{n-1})\in\{o,e\}\times\ldots\times\{o,e\}\backslash\{(e,\dots,e)\}$.
By Definition \ref{DefnC} we obtain
$o=p(e,\dots,e,o)\,\theta\,p(e,\dots,e)=e$. Hence, $o\,\theta\,a$
for all $a\in S$, by compatibility of $\theta$, and we have proved
$\theta=1$.
\end{proof}

\begin{de}\label{ideal}
Let $\bS$ be a semiring. Then $\emptyset\neq I\subseteq S$ is an
ideal of $\bS$ if for all $a,b\in I$ and all $s\in S$ we have
$a+b\in I$ and $as,sa\in I$.
\end{de}

Example: $4\N_0$ is an ideal of $(2\N_0,+,\cdot,0)$.

As in rings, all ideals of a semiring $\bS$ form the lattice
$I(\bS)=(I(S),\cap,+)$. All congruences of $\bS$ form the congruence lattice
$\Con\bS=(\Con\bS,\cap,\vee)$. However, the lattices of ideals and
congruences are not isomorphic. Therefore we deal with certain
congruences that are introduced in \cite[Theorem 8.8]{HW:BookSemirings}. To keep this note self contained we present the proof of Proposition \ref{roi}.

\begin{de}\label{defroi}
Let $\bS$ be a semiring and let $I$ be an ideal of $\bS$. By
$\rho_I$ we denote the smallest congruence of $\bS$ such that $I$ is contained in
the congruence class of $o$.
\end{de}

\begin{pro}\cite[cf. Lemma 8.4, Theorem 8.8]{HW:BookSemirings}\label{roi}
Let $\bS$ be a semiring and let $I$ be an ideal of $\bS$. If we define
$$
a\,\rho\,b\Leftrightarrow(\exists i,j\in I)a+i=b+j
$$
for all $a,b\in S$, then $\rho=\rho_I$.
\end{pro}

\begin{proof}
Obviously, $\rho$ is an equivalence relation on $S$ because $=$ satisfies
reflexivity, symmetry and transitivity. Let us show that $\rho$ is
compatible with $+$ and $\cdot$. Let $a,b\in S$ be such that
$a\,\rho\,b$. By definition of $\rho$ there exist $i,j\in I$
such that $a+i=b+j$. Then for each $c\in S$ we have $c+a+i=c+b+j$.
Hence, $c+a\,\rho\,c+b$ by definition of $\rho$. Therefore
$\rho$ is compatible with $+$. Furthermore, we have
$c(a+i)=c(b+j)$ and $(a+i)c=(b+j)c$. Hence $ca+ci=cb+cj$ and
$ac+ic=bc+jc$ and since $I$ is an ideal we have $ci,cj,ic,jc\in I$
whence $ca\,\rho\,cb$ and $ac\,\rho\,bc$ by definition of
$\rho$. Therefore, $\rho$ is compatible with $\cdot$. Clearly, if $i\in I$ we obtain $i\in[o]_{\rho}$, because $i+o=o+i$. Hence, we have $I\subseteq[o]_{\rho}$. Therefore $\rho_I\subseteq\rho$. If
$(a,b)\in\rho$ and $i,j\in I$ such that $a+i=b+j$ then we have
$a=a+o\,\rho_I\,a+i=b+j\,\rho_I\,b+o=b$. Hence, $(a,b)\in\rho_I$ by
transitivity of $\rho_I$. Therefore, we have proved
$\rho\subseteq\rho_I$ and hence $\rho=\rho_I$.
\end{proof}

We call $\rho_I$ the congruence induced by $I$. The set of all such
congruences of $\bS$ we denote by $ICon\bS$. Note that
$\rho_{\{o\}}=0$ and $\rho_S=1$. Also $I=[o]_{\rho_I}$ if and only if $I$ is $k$-closed, see \cite[Definition 8.5, Theorem 8.8]{HW:BookSemirings}.

\begin{pro}\label{subset}
Let $I$ and $J$ be ideals of a semiring $S$. If $\rho_I=\rho_J$ then $I\subseteq[o]_{\rho_J}$.
\end{pro}

\begin{proof}
By Definition \ref{defroi}.
\end{proof}

Analogously to the product of ideals in rings, we define the product
of finitely many ideals in semirings.

\begin{de}\label{idealproduct}
Let $n\in\N$ and let $I_1,\dots,I_n$ be ideals of a semiring $\bS$.
Then we define
$$
I_1\cdot\ldots\cdot I_n=\{\Sigma_{i=1}^k a^1_i\cdot\ldots\cdot
a^n_i\,|\,k\in\N,a^1_1,\dots,a^1_k\in I_1,\dots,a^n_1,\dots,a^n_k\in
I_n\}.
$$
\end{de}

We observe that the product of finitely many ideals is an ideal and therefore $S^n$ is an ideal of $\bS$ for all $n\in\N$. Similarly as it has been done for rings in \cite[Lemma 3.5]{PM:MASC} we have the following in semirings. 

\begin{pro}\label{commutator ideals}
Let $n\in\N$ and let $I_1,\dots,I_n$ be ideals of a semiring $\bS$.
Then $[I_1,\dots,I_n]$ is an ideal of $\bS$ defined by
$$
[I_1,\dots,I_n]=\Sigma_{\pi\in S_n}I_{\pi(1)}\cdot\ldots\cdot
I_{\pi(n)}.
$$
\end{pro}

\begin{proof}
By definition of ideals and products of ideals in semirings.
\end{proof}

\begin{pro}\label{powersofS}
Let $\bS$ be a semiring and let $m,n,k\in\N$. Then $[S^m,S^n]=S^{m+n}$
and $[\,\underbrace{S,\dots,S}_k\,]=S^k$.
\end{pro}

\begin{proof}
We observe that for every ideal $I$ of $\bS$ we have $I+I=I$ and use
Proposition \ref{commutator ideals}.
\end{proof}

\begin{pro}\label{allabove}
Let $n\in\N$ and let $I_1,\dots,I_n$ be ideals of a semiring $\bS$.
Then for all congruences $\theta$ of $\bS$ we have
$C(\rho_{I_1},\dots,\rho_{I_n};\theta)\Rightarrow\rho_{[I_1,\dots,I_n]}\subseteq\theta$.
\end{pro}

\begin{proof}
Let $C(\rho_{I_1},\dots,\rho_{I_n};\theta)$ and let $i_1\in
I_1,\dots,i_n\in I_n$. Then $o\,\rho_{I_j}\,i_j$ for all
$j\in\{1,\dots,n\}$. For $\pi\in S_n$, let $i_{\pi(1)}\cdot\ldots\cdot i_{\pi(n)}$ be one of the summands of an element from $I_{\pi(1)}\cdot\ldots\cdot I_{\pi(n)}$. Now, we take the polynomial
$p(x_1,\dots,x_n):=x_{\pi(1)}\cdot\ldots\cdot x_{\pi(n)}$. Then we
have $p(x_1,\dots,x_{n-1},i_n)=o\,\theta\,o=p(x_1,\dots,x_{n-1},o)$
for all
$(x_1,\dots,x_{n-1})\in\{o,i_1\}\times\ldots\times\{o,i_{n-1}\}\backslash\{i_1,\dots,i_{n-1}\}$.
Hence we have $p(i_1,\dots,i_n)$ $\theta\,p(i_1,\dots,i_{n-1},o)$.
Therefore $i_{\pi(1)}\cdot\ldots\cdot i_{\pi(n)}=p(i_1,\dots,i_n)\in[o]_{\theta}$. Since every element of
$[I_1,\dots,I_n]$ has finitely many such summands we can easily obtain $[I_1,\dots,I_n]\subseteq[o]_{\theta}$. Now, we let
$a\,\rho_{[I_1,\dots,I_n]}\,b$. Then there exist
$u,v\in[I_1,\dots,I_n]$ such that $a+u=b+v$ by Proposition
\ref{roi}. Hence $o\,\theta\,u$ and $o\,\theta\,v$. Therefore
$a=a+o\,\theta\,a+u=b+v\,\theta\,b+o=b$. Hence $a\,\theta\, b$ by
transitivity of $\theta$. Therefore
$\rho_{[I_1,\dots,I_n]}\subseteq\theta$.
\end{proof}

\begin{cor}\label{oneinequality}
Let $n\in\N$ and let $I_1,\dots,I_n$ be ideals of a semiring $\bS$.
Then we have
$\rho_{[I_1,\dots,I_n]}\leq[\rho_{I_1},\dots,\rho_{I_n}]$.
\end{cor}

\begin{proof}
The inequality follows from Definition \ref{DefnHC} and Proposition \ref{allabove}.
\end{proof}

Let us note that the opposite inequality in Corollary \ref{oneinequality} is not always true. We take the semiring $\mathbf B=(\{\top,\bot\},\vee,\cdot)$ such that the result of $\cdot$ is always $\bot$. The semiring $\mathbf B$ has only two congruences $0$ and $1$. Since the disjunction $\vee$ is a semilattice operation we know that the commutator of two congruences is their intersection by \cite[Exercise 5]{McMc}. Therefore, $[\rho_{\{\top,\bot\}},\rho_{\{\top,\bot\}}]=[1,1]=1\neq0=\rho_{\{\bot\}}=\rho_{\{\top,\bot\}^2}=\rho_{[\{\top,\bot\},\{\top,\bot\}]}$.

\section{Properties Defined by Commutators}

As usual, for each polynomial function of a semiring $\bS$ there is
a polynomial word (polynomial term), a term with constants that
induces the given polynomial function. Without confusion we call both
a polynomial word and a polynomial function a polynomial. In a semiring
$\bS=(S,+,\cdot,o)$ polynomials of $(S,\cdot)$ we call
\emph{monomials}. The following proposition we are going to use in the sequel without explicit reference.

\begin{pro}\label{sumofmonomials}
In semirings each polynomial can be written as a sum of monomials.
\end{pro}

\begin{proof}
One can prove that a polynomial with the polynomial word
$p(x_1,\dots,x_n)$ has the polynomial word in the form of the sum of
monomials, by induction of the number of occurrences of $+$ and
$\cdot$ in $p(x_1,\dots,x_n)$.
\end{proof}

\begin{pro}\label{withcancelativity}
Let $I$ and $J$ be ideals of a semiring $\bS$. Then
$C(\rho_I,\rho_J;\rho_{[I,J]})$ if and only if we have
$a+c\,\rho_{[I,J]}\,a+d\Rightarrow b+c\,\rho_{[I ,J]}\,b+d$ for all
$a,b,c,d\in S$ such that $a\,\rho_I\,b$ and $c\,\rho_J\,d$.
\end{pro}

\begin{proof} ($\Rightarrow$) We assume $C(\rho_I,\rho_J;\rho_{[I,J]})$ and apply Definition \ref{DefnC} for the polynomial $p(x,y)=x+y$ to obtain the statement. ($\Leftarrow$) We prove
$C(\rho_I,\rho_J;\rho_{[I,J]})$ using Proposition \ref{DefnC(a,a,a,b;d)}. We let 
$p(a,\oc)\,\rho_{[I,J]}\,p(a,\od)$ for an ($n+1$)-ary polynomial
$p(x,\oy)$ and $a\,\rho_I\,b$ and $\oc\,\rho_J\,\od$. By Proposition
\ref{sumofmonomials}, we know that there exist monomials $s(x)$ and
$r_1(x,\oy),\dots,r_k(x,\oy)$ such that
$p(x,\oy)=s(x)+r_1(x,\oy)+\ldots+r_k(x,\oy)$, where each of
$r_1(x,\oy),\dots,r_k(x,\oy)$ may depend on $x$. We prove the
statement $p(b,\oc)\,\rho_{[I,J]}\,p(b,\od)$ by induction on the
number $m\in\N_0$ of occurrences of $x$ in the polynomial word of
$r_1(x,\oy)+\ldots+r_k(x,\oy)$. If $m=0$ then we denote
$r_1(a,\oc)+\ldots+r_k(a,\oc)=r_1(b,\oc)+\ldots+r_k(b,\oc)$ by
$r(\oc)$ and
$r_1(a,\od)+\ldots+r_k(a,\od)=r_1(b,\od)+\ldots+r_k(b,\od)$ by
$r(\od)$ and the statement is $
s(a)+r(\oc)\,\rho_{[I,J]}\,s(a)+r(\od)\Rightarrow
s(b)+r(\oc)\,\rho_{[I,J]}\,s(b)+r(\od). $ It is true by the
assumption, because $s(a)\,\rho_I\,s(b)$ and
$r(\oc)\,\rho_J\,r(\od)$. Let $m\geq 1$ and let us assume that the variable $x$ appears in the monomial word for $r_1(x,\oy)$, otherwise we rename the indexes of the monomials. We
 denote
$r_2(x,\oy)+\ldots+r_k(x,\oy)$ by $r(x,\oy)$ and analyze the following cases.

(1) $r_1(x,\oy)=xq(x,\oy)$, where $q(x,\oy)$ is a submonomial of $r_1(x,\oy)$. Let
$c:=q(a,\oc)\,\rho_J\,q(a,\od)=:d$. Then, there exist $i,j\in I$
and $k,\ell\in J$ such that $a+i=b+j$ and $c+k=d+\ell$. Further, we have 
assumed $p(a,\oc)\,\rho_{[I,J]}\,p(a,\od)$. This is
$s(a)+ac+r(a,\oc)\,\rho_{[I,J]}\,s(a)+ad+r(a,\od)$. Then, we
have
$$
s(a)+ac+ic+id+r(a,\oc)\,\rho_{[I,J]}\,s(a)+ad+ic+id+r(a,\od).
$$
This yields
$$
s(a)+(a+i)c+id+r(a,\oc)\,\rho_{[I,J]}\,s(a)+(a+i)d+ic+r(a,\od).
$$
Therefore,
$$
s(a)+(b+j)c+id+r(a,\oc)\,\rho_{[I,J]}\,s(a)+(b+j)d+ic+r(a,\od).
$$
Hence,
$$
s(a)+bc+jc+id+r(a,\oc)\,\rho_{[I,J]}\,s(a)+bd+jd+ic+r(a,\od).
$$
Thus,
$$
s(a)+bc+jc+jk+id+i\ell+j\ell+ik+r(a,\oc)\,
$$
$$
\rho_{[I,J]}\,s(a)+bd+jd+j\ell+ic+ik+jk+i\ell+r(a,\od).
$$
Whence,
$$
s(a)+bc+j(c+k)+i(d+\ell)+j\ell+ik+r(a,\oc)\,
$$
$$
\rho_{[I,J]}\,s(a)+bd+j(d+\ell)+i(c+k)+jk+i\ell+r(a,\od).
$$
Since $d+\ell=c+k$ we have $j(c+k)=j(d+\ell)$ and
$i(d+\ell)=i(c+k)$. Hence,
$j(c+k)+i(d+\ell)=j(d+\ell)+i(c+k)\,\rho_I\,o$. By the
assumption we obtain
$$
s(a)+bc+j\ell+ik+r(a,\oc)\,\rho_{[I,J]}\,s(a)+bd+jk+i\ell+r(a,\od),
$$
because $j\ell+ik\in J$ and $jk+i\ell\in J$ and $s(a)+bc+r(a,\oc)\,\rho_{J}\,s(a)+bd+r(a,\od)$ and therefore
$s(a)+bc+j\ell+ik+r(a,\oc)\,\rho_{J}\,s(a)+bd+jk+i\ell+r(a,\od)$.
Since $j\ell+ik\,\rho_{[I,J]}\,jk+i\ell\,\rho_{[I,J]}\,o$
then we have
$s(a)+bc+r(a,\oc)\,\rho_{[I,J]}\,s(a)+bd+r(a,\od)$ by transitivity of $\rho_{[I,J]}$. Hence we
obtain
$s(a)+bq(a,\oc)+r(a,\oc)\,\rho_{[I,J]}\,s(a)+bq(a,\od)+r(a,\od)$.
The polynomial $bq(x,\oy)+r(x,\oy)$ has $m-1$ occurrences of
$x$ in its polynomial word and by the induction hypothesis we have
$s(b)+bq(b,\oc)+r(b,\oc)\,\rho_{[I,J]}\,s(b)+bq(b,\od)+r(b,\od)$.
Therefore, $p(b,\oc)\,\rho_{[I,J]}\,p(b,\od)$.

(2)  $r_1(x,\oy)=t(\oy)xq(x,\oy)$, where $t(\oy)$ and $q(x,\oy)$ are submonomials of $r_1(x,\oy)$ such that $t(\oy)$ does not depend on $x$. If we denote $t(\oc)$ by $c'$ and $t(\od)$ by $d'$, we have $c'\,\rho_J\,d'$ and we let
$c:=q(a,\oc)\,\rho_J\,q(a,\od)=:d$. Then, there exist $i,j\in I$
and $k,k',\ell,\ell'\in J$ such that $a+i=b+j$, $c+k=d+\ell$ and $c'+k'=d'+\ell'$. Further, we have
assumed $p(a,\oc)\,\rho_{[I,J]}\,p(a,\od)$. This is
$s(a)+c'ac+r(a,\oc)\,\rho_{[I,J]}\,s(a)+d'ad+r(a,\od)$. Then, we
have
$$
s(a)+c'ac+c'ic+d'id+r(a,\oc)\,\rho_{[I,J]}\,s(a)+d'ad+c'ic+d'id+r(a,\od).
$$
This yields
$$
s(a)+c'(a+i)c+d'id+r(a,\oc)\,\rho_{[I,J]}\,s(a)+d'(a+i)d+c'ic+r(a,\od).
$$
Therefore,
$$
s(a)+c'(b+j)c+d'id+r(a,\oc)\,\rho_{[I,J]}\,s(a)+d'(b+j)d+c'ic+r(a,\od).
$$
Hence,
$$
s(a)+c'bc+c'jc+d'id+r(a,\oc)\,\rho_{[I,J]}\,s(a)+d'bd+d'jd+c'ic+r(a,\od).
$$
Thus,
$$
s(a)+c'bc+c'jc+c'jk+d'id+d'i\ell+d'j\ell+c'ik+r(a,\oc)\,
$$
$$
\rho_{[I,J]}\,s(a)+d'bd+d'jd+d'j\ell+c'ic+c'ik+c'jk+d'i\ell+r(a,\od).
$$
Whence,
$$
s(a)+c'bc+c'j(c+k)+d'i(d+\ell)+d'j\ell+c'ik+r(a,\oc)\,
$$
$$
\rho_{[I,J]}\,s(a)+d'bd+d'j(d+\ell)+c'i(c+k)+c'jk+d'i\ell+r(a,\od).
$$
Hence we obtain
$$
s(a)+c'bc+c'j(c+k)+k'j(c+k)+d'i(d+\ell)+\ell'i(d+\ell)+d'j\ell+c'ik+\ell'j(d+\ell)+k'i(c+k)+r(a,\oc)\,
$$
$$
\rho_{[I,J]}
$$
$$
s(a)+d'bd+d'j(d+\ell)+\ell'j(d+\ell)+c'i(c+k)+k'i(c+k)+c'jk+d'i\ell+k'j(c+k)+\ell'i(d+\ell)+r(a,\od).
$$
Whence,
$$
s(a)+c'bc+(c'+k')j(c+k)+(d'+\ell')i(d+\ell)+d'j\ell+c'ik+\ell'j(d+\ell)+k'i(c+k)+r(a,\oc)\,
$$
$$
\rho_{[I,J]}\,s(a)+d'bd+(d'+\ell')j(d+\ell)+(c'+k')i(c+k)+c'jk+d'i\ell+k'j(c+k)+\ell'i(d+\ell)+r(a,\od).
$$
Since $c'+k'=d'+\ell'$ and  $c+k=d+\ell$ we have $(c'+k')j(c+k)=(d'+\ell')j(d+\ell)$ and
$(c'+k')i(c+k)=(d'+\ell')i(d+\ell)$. Hence,
$(c'+k')j(c+k)+(d'+\ell')i(d+\ell)=(d'+\ell')j(d+\ell)+(c'+k')i(c+k)\,\rho_I\,o$. By the
assumption we obtain
$$
s(a)+c'bc+d'j\ell+c'ik+\ell'j(d+\ell)+k'i(c+k)+r(a,\oc)\,
$$
$$
\rho_{[I,J]}\,s(a)+d'bd+c'jk+d'i\ell+k'j(c+k)+\ell'i(d+\ell)+r(a,\od),
$$
because $d'j\ell+c'ik+\ell'j(d+\ell)+k'i(c+k)\in J$, $c'jk+d'i\ell+k'j(c+k)+\ell'i(d+\ell)\in J$ and $s(a)+c'bc+r(a,\oc)\,\rho_{J}\,s(a)+d'bd+r(a,\od)$ and therefore
the both sides of the last relation are in $\rho_{J}$.
Since $d'j\ell+c'ik+\ell'j(d+\ell)+k'i(c+k)\,\rho_{[I,J]}\,c'jk+d'i\ell+k'j(c+k)+\ell'i(d+\ell)\,\rho_{[I,J]}\,o$
then we have
$s(a)+c'bc+r(a,\oc)\,\rho_{[I,J]}\,s(a)+d'bd+r(a,\od)$ by transitivity of $\rho_{[I,J]}$. Hence we
obtain
$s(a)+t(\oc)bq(a,\oc)+r(a,\oc)\,\rho_{[I,J]}\,s(a)+t(\od)bq(a,\od)+r(a,\od)$.
The polynomial $t(\oy)bq(x,\oy)+r(x,\oy)$ has $m-1$ occurrences of
$x$ in its polynomial word and by the induction hypothesis we have
$s(b)+t(\oc)bq(b,\oc)+r(b,\oc)\,\rho_{[I,J]}\,s(b)+t(\od)bq(b,\od)+r(b,\od)$.
Therefore, $p(b,\oc)\,\rho_{[I,J]}\,p(b,\od)$.

(3) $r_1(x,\oy)=q(x,\oy)x$ is dual to the case (1).
\end{proof}

\begin{cor}\label{allworks}
In an additively cancellative semiring $\bS=(S,+,\cdot,o)$ we have $[\rho_I,\rho_J]=\rho_{[I,J]}$, for all ideals $I$ and $J$ of $\bS$ .
\end{cor}

\begin{proof}
Let $\bS$ be an additively cancellative semiring. We take $a,b,c,d\in S$ such that $a\,\rho_I\,b$ and $c\,\rho_J\,d$ and assume $a+c\,\rho_{[I,J]}\,a+d$. By Proposition \ref{roi} we know that there exist $i,j\in[I,J]$ such that $a+c+i=a+d+j$. Hence, $c+i=d+j$ by cancellativity for $+$. Therefore, $b+c+i=b+d+j$ and we obtain $b+c\,\rho_{[I,J]}\,b+d$ again by Proposition \ref{roi}. Now, we have $C(\rho_I,\rho_J;\rho_{[I,J]})$ by Proposition \ref{withcancelativity}. Hence, $[\rho_I,\rho_J]\leq\rho_{[I,J]}$ by Definition \ref{DefnHC}. The opposite inequality has been proved in Corollary \ref{oneinequality}.
\end{proof}

\begin{pro}\label{teoremasolvable}
Let $k\in\N$. An additively cancellative semiring $\bS=(S,+,\cdot,o)$ is $k$-solvable if and only if $S^{2^k}=\{o\}$.
\end{pro}

\begin{proof}
  One can prove that $[1,1]^{(k)}=[\rho_S,\rho_S]^{(k)}=\rho_{S^{2^k}}$ by induction on $k$, using Definition \ref{DefNilpSolv} and Corollary \ref{allworks}.
\end{proof}

\begin{lem}\label{fornilpotency}
If a commutative monoid $(S,+,o)$ is nilpotent (supernilpotent), then it is cancellative.
\end{lem}

\begin{proof}
First, we assume that $(S,+,o)$ is $k$-nilpotent for $k\in\N$. Furthermore, we take $a,c,d\in S$ such that $a+c=a+d$. Then $a+c\,[1,1]\,a+d$. Using Proposition \ref{DefnC(a,a,a,b;d)} for the polynomial $p(x,y)=x+y$ we obtain $c=o+c\,[1,1]\,o+d=d$. Similarly, from the assumption $a+c=a+d$, we have $a+c\,[1,[1,1]]\,a+d$ and thus, $c\,[1,[1,1]]\,d$. After $k$ steps we obtain $c\,(1,1]^{(k)}\,d$ and hence $c=d$.

Now, we assume that $(S,+,o)$ is $n$-supernilpotent for $n\in\N$. Furthermore, we take $a,c,d\in S$ such that  $a+c=a+d$. Then $ka+c=ka+d$ and hence $ka+c\,[\,\underbrace{1,\dots,1}_{n+1}\,]\,ka+d$, for all $k\in\{1,\dots,n\}$. Using Definitions \ref{DefnC} and \ref{DefnHC} for the polynomial $p(x_1,\dots,x_{n},y)=x_1+\ldots+x_{n}+y$ we obtain $c=o+c\,[\,\underbrace{1,\dots,1}_{n+1}\,]\,o+d=d$ and hence $c=d$.
\end{proof}

 For a semiring $(S,+,\cdot,o)$, $n\in\N$, tuples $\oa_1,\dots,\oa_n,\ob_1,\dots,\ob_n$ from $S$ and tuples of variables $\ox_1,\dots,\ox_n$ such that $\oa_i,\ob_i$ and $\ox_i$ have the same length for all $i\in\{1,\dots,n\}$, we denote the set of valuations $v:\{(\ox_1,\dots,\ox_n)\}\to\{\oa_1,\ob_1\}\times\dots\times\{\oa_n,\ob_n\}$ by $V^n_{\oa,\ob}$ and $|\{i\in\{1,\dots,n\}\,:\,v(\ox_i)=\ob_i\}|$ by $|v|_b$. According to the Proposition \ref{sumofmonomials}, for each polynomial we have $p(\ox_1,\dots,\ox_n)=\Sigma^s_{\ell=1}m_{\ell}(\ox_{i_1},\dots,\ox_{i_k})$, where $s\in\N$ and $m(\ox_{i_1},\dots,\ox_{i_k})$ are monomials that contain at least one component of each of the tuples of variables $\ox_{i_1},\dots,\ox_{i_k}$ at least once. 
 
 \begin{lem}\label{monomial}
Let $k,n\in\N$ be such that  $n\geq 2$ and $k<n$. Let $(S,+,\cdot,o)$ be a semiring, and let $\ox_1,\dots,\ox_n$ be tuples of variables and $\oa_1,\dots,\oa_n,\ob_1,\dots,\ob_n$ tuples from $S$ such that $\oa_i,\ob_i$ and $\ox_i$ have the same length for all $i\in\{1,\dots,n\}$. If $1\leq i_1<\dots<i_k\leq n$ and $m(\ox_{i_1},\dots,\ox_{i_k})$ is a monomial that contains at least one component of each of the tuples of variables $\ox_{i_1},\dots,\ox_{i_k}$ at least once, then we have 
 $$
 \Sigma_{v\in V^n_{\oa,\ob},|v|_b\equiv_2 0}m^v(\ox_{i_1},\dots,\ox_{i_k})= \Sigma_{v\in V^n_{\oa,\ob},|v|_b\equiv_2 1}m^v(\ox_{i_1},\dots,\ox_{i_k}).
 $$
\end{lem}

\begin{proof}
We prove the statement by induction on $n$. For $n=2$ we have $k=1$. There are only 4 valuations: $v_1=\left(\begin{array}{cc}\ox_1&\ox_2\\ \oa_1&\oa_2\end{array}\right)$, $v_2=\left(\begin{array}{cc}\ox_1&\ox_2\\ \ob_1&\oa_2\end{array}\right)$, $v_3=\left(\begin{array}{cc}\ox_1&\ox_2\\ \oa_1&\ob_2\end{array}\right)$ and $v_4=\left(\begin{array}{cc}\ox_1&\ox_2\\ \ob_1&\ob_2\end{array}\right)$. Let $i_1=1$, the other case is analogous. Now, we have $m^{v_1}(\ox_1)+m^{v_4}(\ox_1)=m(\oa_1)+m(\ob_1)=m^{v_3}(\ox_1)+m^{v_2}(\ox_1)$, because $|v_1|_b=0$ and $|v_4|_b=2$ are even and  $|v_2|_b=|v_3|_b=1$ are odd. We let $n\geq 3$. First, we assume that $k<n-1$. If $\ox_n\not\in\{\ox_{i_1},\dots,\ox_{i_k}\}$ the statement follows by the induction hypothesis. Otherwise, we rename the variables. Let us now assume $k=n-1$ and $\{\ox_{i_1},\dots,\ox_{i_k}\}=\{\ox_1,\dots,\ox_{n-1}\}$ without loss of generality. For each valuation $v$ there exist the valuation $v'$ such that $v'(\ox_i)=v(\ox_i)$ for all $i\in\{1,\dots,n-1\}$ and $v'(\ox_n)\neq v(\ox_n)$. Let us observe that $m^v\tup{\ox}{1}{n-1}=m^{v'}\tup{\ox}{1}{n-1}$, but $|v|_b\not\equiv_2 |v'|_b$. Therefore $m^v\tup{\ox}{1}{n-1}$ and $m^{v'}\tup{\ox}{1}{n-1}$ appear on the opposite sides of the equality in the statement. Hence, the equality is true. This finishes the proof of the induction step.
\end{proof}

\begin{cor}\label{polynomial}
Let $k,n\in\N$ be such that $n\geq 2$ and $k<n$. Let $p(\ox_1,\dots,\ox_n)$ be a polynomial over a semiring $(S,+,\cdot,o)$ such that $p(\ox_1,\dots,\ox_n)=\Sigma^s_{\ell=1} m_{\ell}(\ox_{i_1},\dots,\ox_{i_k})$, where $s\in\N$ and $m_{\ell}(\ox_{i_1},\dots,\ox_{i_k})$ are monomials that contain at least one component of each of the variables $\ox_{i_1},\dots,\ox_{i_k}$ at least once. If $\oa_1,\dots,\oa_n,\ob_1,\dots,\ob_n$ are tuples from $S$ such that $\oa_i,\ob_i$ and $\ox_i$ have the same length for all $i\in\{1,\dots,n\}$ then, we have
$$
 \Sigma_{v\in V^n_{\oa,\ob},|v|_b\equiv_2 0}p^v(\ox_1,\dots,\ox_n)= \Sigma_{v\in V^n_{\oa,\ob},|v|_b\equiv_2 1}p^v(\ox_1,\dots,\ox_n).
 $$
\end{cor}

\begin{proof}
By commutativity and associatvity for $+$ and Lemma \ref{monomial} we obtain:

$\Sigma_{v\in V^n_{\oa,\ob},|v|_b\equiv_2 0}p^v(\ox_1,\dots,\ox_n)=\Sigma_{v\in V^n_{\oa,\ob},|v|_b\equiv_2 0}(\Sigma^s_{\ell=1} m^v_{\ell}(\ox_{i_1},\dots,\ox_{i_k}))=$ 

$\Sigma^s_{\ell=1}(\Sigma_{v\in V^n_{\oa,\ob},|v|_b\equiv_2 0} m^v_{\ell}(\ox_{i_1},\dots,\ox_{i_k}))=\Sigma^s_{\ell=1}(\Sigma_{v\in V^n_{\oa,\ob},|v|_b\equiv_2 1} m^v_{\ell}(\ox_{i_1},\dots,\ox_{i_k}))=$ 

$\Sigma_{v\in V^n_{\oa,\ob},|v|_b\equiv_2 1}(\Sigma^s_{\ell=1} m^v_{\ell}(\ox_{i_1},\dots,\ox_{i_k}))=\Sigma_{v\in V^n_{\oa,\ob},|v|_b\equiv_2 1}p^v(\ox_1,\dots,\ox_n)$.
\end{proof}

Now, we prove the main result, Theorem \ref{teoremanilpotent}.

\begin{proof}
$(1)\Rightarrow(2)$ We have that the commutative monoid $(S,+,o)$ is also $n$-nilpotent by Proposition \ref{forreduct}, because $\{+,o\}\subseteq\{+,\cdot,o\}$. Then, using Lemma \ref{fornilpotency} we obtain that $(S,+,o)$ is cancellative and therefore, $\bS$ is additively cancellative. Thus, one can prove that $(1,1]^{(k)}=\rho_{S^{k+1}}$ by induction on $k$, using Definition \ref{DefNilpSolv} and Corollary \ref{allworks}. Hence $\rho_{S^{n+1}}=(1,1]^{(n)}=0=\rho_{\{o\}}$ and therefore $S^{n+1}=\{o\}$ by Proposition \ref{subset}. 

$(2)\Rightarrow(3)$ Since $S^{n+1}=\{o\}$, then for every polynomial $p(\ox_1,\dots,\ox_n,\oy)=q\tup{\ox}{1}{n}+\Sigma^s_{\ell=1}m_{\ell}(\ox_{i_1},\dots,\ox_{i_k},\oy)$, where $m_{\ell}(\ox_{i_1},\dots,\ox_{i_k},\oy)$ are monomials that contain at least one component of each of the variables $\ox_{i_1},\dots,\ox_{i_k},\oy$ at least once, we have $k<n$. Let $n$ be an even number, the opposite case is analogous. We let $\oa_1,\dots,\oa_n,\ob_1,\dots,\ob_n,\oc,\od$ tuples in $S$ be such that $p^v(\ox_1,\dots,\ox_n,\oc)=p^v(\ox_1,\dots,\ox_n,\od)$ for all valuations $v\in V^n_{\oa,\ob}$ such that there exist a $j\in\{1,\dots,n\}$ with  $v(\ox_j)=\oa_j$. Hence, if we define $w(\ox_i)=\ob_i$ for all $i\in\{1,\dots,n\}$, we obtain 
$$\Sigma_{v\neq w,|v|_b\equiv_2 0}p^v(\ox_1,\dots,\ox_n,\oc)=\Sigma_{v\neq w,|v|_b\equiv_2 0}p^v(\ox_1,\dots,\ox_n,\od)$$ and 
$$\Sigma_{v\in V^n_{\oa,\ob},|v|_b\equiv_2 1}p^v(\ox_1,\dots,\ox_n,\oc)=\Sigma_{v\in V^n_{\oa,\ob},|v|_b\equiv_2 1}p^v(\ox_1,\dots,\ox_n,\od).$$  Now, using  Co\-ro\-llary \ref{polynomial} we have:

$p(\ob_1,\dots,\ob_n,\oc)+\Sigma_{v\neq w,|v|_b\equiv_2 0}p^v(\ox_1,\dots,\ox_n,\oc)$ 

$=\Sigma_{v\in V^n_{\oa,\ob},|v|_b\equiv_2 0}p^v(\ox_1,\dots,\ox_n,\oc)=\Sigma_{v\in V^n_{\oa,\ob},|v|_b\equiv_2 1}p^v(\ox_1,\dots,\ox_n,\oc)$ 

$=\Sigma_{v\in V^n_{\oa,\ob},|v|_b\equiv_2 1}p^v(\ox_1,\dots,\ox_n,\od)=\Sigma_{v\in V^n_{\oa,\ob},|v|_b\equiv_2 0}p^v(\ox_1,\dots,\ox_n,\od)$ 

$=p(\ob_1,\dots,\ob_n,\od)+\Sigma_{v\neq w,|v|_b\equiv_2 0}p^v(\ox_1,\dots,\ox_n,\od)$.

Therefore, $p(\ob_1,\dots,\ob_n,\oc)=p(\ob_1,\dots,\ob_n,\od)$ by cancellativity of $+$. Hence, we obtain $C(\,\underbrace{1,\dots,1}_{n+1};0)$ and the semiring is $n$-supernilpotent by definition.

$(3)\Rightarrow(1)$ Let $\bS$ be $n$-supernilpotent. Then, by Proposition \ref{forreduct}, the commutative monoid $(S,+,o)$ is also $n$-supernilpotent, because $\{+,o\}\subseteq\{+,\cdot,o\}$. Therefore, the monoid $(S,+,o)$ is cancellative by Lemma \ref{fornilpotency} and hence $\bS$ is additively cancellative. Also, we have $0=[\,\underbrace{1,\dots,1}_{n+1}\,]=[\,\underbrace{\rho_S,\dots,\rho_S}_{n+1}\,]\geq\rho_{[\,\underbrace{S,\dots,S}_{n+1}\,]}=\rho_{S^{n+1}}$, by Corollary \ref{oneinequality} and Proposition \ref{powersofS}. Since $\bS$ is additively cancellative we can prove $(1,1]^{(k)}=(\rho_S,\rho_S]^{(k)}=\rho_{S^{k+1}}$ by induction on $k$, using Definition \ref{DefNilpSolv} and Corollary \ref{allworks}. Hence, $\bS$ is $n$-nilpotent because $\rho_{S^{n+1}}=0$.
\end{proof}

\begin{cor}\label{finitenilpotent}
Let $k\in\N$. A finite semiring $\bS=(S,+,\cdot,o)$ is $k$-supernilpotent ($k$-nilpotent) if and only if $\bS$ is a $k$-nilpotent ring.
\end{cor}

\begin{proof}
$(\Rightarrow)$ Using Theorem \ref{teoremanilpotent} and the well known fact that a finite cancellative monoid is a group we obtain that $\bS=(S,+,\cdot,o)$ is a finite ring with $S^{k+1}=\{o\}$. Since rings are additively cancellative semirings, we obtain that $\bS=(S,+,\cdot,o)$ is a finite $k$-supernilpotent ($k$-nilpotent) ring, by Theorem \ref{teoremanilpotent}. $(\Leftarrow)$ Obvious, since every ring is an additively cancellative semiring and therefore it is $k$-supernilpotent if and only if it is $k$-nilpotent by Theorem \ref{teoremanilpotent}.
\end{proof}

In the following statement we denote the ideal $\{s_1\cdot\ldots\cdot s_{n+1}\,|\,s_1,\dots,s_{n+1}\in S\}$ of the semigroup $\bS$ by $S^{n+1}$ as in \cite{RM:cswz}. If $\bS$ is a semiring we observe that the condition $S^{n+1}=\{o\}$ in this sense and in the sense of the Definition \ref{idealproduct} is the same.

\begin{pro}\rm{(cf. \cite[Proposition 3.5]{RM:cswz}).}\label{jelininimojrad}
Let $n\in\N$. A semigroup with zero $(S,\cdot,o)$ is $n$-supernilpotent ($n$-nilpotent) if and only if $S^{n+1}=\{o\}$.
\end{pro}

\begin{pro}\label{bothnilpotent}
Let $k\in\N$. A semiring $\bS=(S,+,\cdot,o)$ is $k$-supernilpotent ($k$-nilpotent) if and only if the monoid $(S,+,o)$ and the semigroup with zero $(S,\cdot,o)$ are $k$-supernilpotent ($k$-nilpotent).
\end{pro}

\begin{proof}
  $(\Rightarrow)$ By the assumption we obtain that commutative monoid $(S,+,o)$ and the semigroup with zero $(S,\cdot,o)$ are also $k$-supernilpotent ($k$-nilpotent) by Proposition \ref{forreduct}, because $\{+,o\},\{\cdot,o\}\subseteq\{+,\cdot,o\}$. $(\Leftarrow)$ From supernilpotency (nilpotency) of $(S,+,o)$ we obtain that $\bS$ is additively cancellative by Lemma \ref{fornilpotency}, and from $k$-supernilpotency ($k$-nilpotency) of $(S,\cdot,o)$ we obtain that $S^{k+1}=\{o\}$ by Proposition \ref{jelininimojrad} and therefore $S^{k+1}=\{o\}$ in the sense of Definition \ref{idealproduct}. The statement follows by Theorem \ref{teoremanilpotent}.
\end{proof}

\begin{pro}\label{teoremabelian}
A semiring $\bS=(S,+,\cdot,o)$ is abelian if and only if $\bS$ is additively
cancellative and the multiplication $\cdot$ is the zero multiplication.
\end{pro}

\begin{proof}
This is a corollary of Theorem \ref{teoremanilpotent} using the fact that abelian means $1$-nilpotent. 
\end{proof}

\begin{cor}\label{finiteabelian}
A finite semiring is abelian if and only if it is a zero-ring.
\end{cor}

\begin{proof}
The statement follows from Proposition \ref{teoremabelian} using the fact that finite cancellative monoid is a group.
\end{proof}

\section{Acknowledgements}

The authors thank to P. Mayr for useful comments during preparation of the manuscript.

\bibliographystyle{plain}

\end{document}